\documentclass[11pt]{amsart}
\usepackage{amsfonts,amssymb,amscd,amsmath,enumerate,verbatim,calc,latexsym,pstcol,pst-plot,pst-3d,}

\def\NZQ{\mathbb}               
\def\NN{{\NZQ N}}

\def\ZZ{{\NZQ Z}}
\def\RR{{\NZQ R}}
\def\CC{{\NZQ C}}

\def\frk{\mathfrak}               

\def\Phi{{\frk N}}
\def\opn#1#2{\def#1{\operatorname{#2}}} 
\opn\chara{char} \opn\length{\ell} \opn\pd{pd} \opn\rk{rk}
\opn\projdim{proj\,dim} \opn\injdim{inj\,dim} \opn\rank{rank}
\opn\depth{depth} \opn\grade{grade} \opn\height{height}
\opn\embdim{emb\,dim} \opn\codim{codim}

\opn\Tr{Tr} \opn\bigrank{big\,rank}
\opn\superheight{superheight}\opn\lcm{lcm}
\opn\trdeg{tr\,deg}
\opn\reg{reg} \opn\lreg{lreg} \opn\ini{in} \opn\lpd{lpd}
\opn\size{size}\opn{\mult}{mult}
\opn\div{div} \opn\Div{Div} \opn\cl{cl} \opn\Cl{Cl}
\opn\Spec{Spec} \opn\Supp{Supp} \opn\supp{supp} \opn\Sing{Sing}
\opn\Ass{Ass} \opn\Min{Min}
\opn\Ann{Ann} \opn\Rad{Rad} \opn\Soc{Soc}
\opn\Syz{Syz} \opn\Im{Im} \opn\Ker{Ker} \opn\Coker{Coker}
\opn\Am{Am} \opn\Hom{Hom} \opn\Tor{Tor} \opn\Ext{Ext}
\opn\End{End} \opn\Aut{Aut} \opn\id{id} \opn\ini{in}

\opn\nat{nat}
\opn\pff{pf}
\opn\Pf{Pf} \opn\GL{GL} \opn\SL{SL} \opn\mod{mod} \opn\ord{ord}
\opn\Gin{Gin}
\opn\Hilb{Hilb}\opn\adeg{adeg}\opn\std{std}\opn\ip{infpt}
\opn\Pol{Pol}
\opn\sat{sat}
\opn\Var{Var}
\opn\Gen{Gen}

\opn\aff{aff} \opn\con{conv} \opn\relint{relint} \opn\st{st}
\opn\lk{lk} \opn\cn{cn} \opn\core{core} \opn\vol{vol}
\opn\link{link} \opn\star{star}
\opn\gr{gr}

\def\pot#1#2{#1[\kern-0.28ex[#2]\kern-0.28ex]}

\opn\dirlim{\underrightarrow{\lim}}
\opn\inivlim{\underleftarrow{\lim}}
\let\union=\cup
\let\sect=\cap

\let\iso=\cong
\let\Union=\bigcup
\let\Sect=\bigcap

\let\to=\rightarrow

\def\Implies{\ifmmode\Longrightarrow \else
        \unskip${}\Longrightarrow{}$\ignorespaces\fi}
\def\implies{\ifmmode\Rightarrow \else
        \unskip${}\Rightarrow{}$\ignorespaces\fi}
\def\iff{\ifmmode\Longleftrightarrow \else
        \unskip${}\Longleftrightarrow{}$\ignorespaces\fi}

\let\:=\colon
\newtheorem{Theorem}{Theorem}[section]
\newtheorem{Lemma}[Theorem]{Lemma}
\newtheorem{Corollary}[Theorem]{Corollary}
\newtheorem{Proposition}[Theorem]{Proposition}
\newtheorem{Remark}[Theorem]{Remark}

\newtheorem{Example}[Theorem]{Example}
\newtheorem{Examples}[Theorem]{Examples}
\newtheorem{Definition}[Theorem]{Definition}

\let\epsilon\varepsilon
\let\phi=\varphi
\let\kappa=\varkappa
\def\qed{\ifhmode\textqed\fi
      \ifmmode\ifinner\quad\qedsymbol\else\dispqed\fi\fi}
\def\textqed{\unskip\nobreak\penalty50
       \hskip2em\hbox{}\nobreak\hfil\qedsymbol
       \parfillskip=0pt \finalhyphendemerits=0}
\def\dispqed{\rlap{\qquad\qedsymbol}}

\opn\dis{dis}
\def\pnt{{\raise0.5mm\hbox{\large\bf.}}}

\opn\Lex{Lex}

\newcommand\Perp{\protect\mathpalette{\protect\independenT}{\perp}}
\def\independenT#1#2{\mathrel{\rlap{$#1#2$}\mkern2mu{#1#2}}}
\newcommand{\ind}[2]{\left.#1 \Perp #2 \inD}
\newcommand{\inD}[1][\relax]{\def\argone{#1}\def\temprelax{\relax}
  \ifx\argone\temprelax\right.\else\,\middle|#1\right.{}\fi}
\newcommand{\CI}{CI-}

\newif\ifbinary
\binarytrue

\begin{document}

\title{Binomial edge ideals and conditional independence statements}

\author{J\"urgen Herzog, Takayuki Hibi, Freyja Hreinsd\'{o}ttir, Thomas Kahle, and Johannes Rauh}

\address{J\"urgen Herzog, Fachbereich Mathematik, Universit\"at Duisburg-Essen, Campus Essen, 45117
Essen, Germany} \email{juergen.herzog@uni-essen.de}

\address{Takayuki Hibi, Department of Pure and Applied Mathematics, Graduate School of Information Science and Technology,
Osaka University, Toyonaka, Osaka 560-0043, Japan}
\email{hibi@math.sci.osaka-u.ac.jp}

\address{Freyja Hreinsd\'{o}ttir,   School of Education, University of Iceland, Stakkahlid, 105 Reykjavik, Iceland}
\email{freyjah@hi.is}

\address{Thomas Kahle and Johannes Rauh, MPI for Mathematics in the Sciences, 04103, Leipzig, Germany}
\email{\{kahle,rauh\}@mis.mpg.de}

\begin{abstract}
  We introduce binomial edge ideals attached to a simple graph $G$ and study their algebraic properties.  We
  characterize those graphs for which the quadratic generators form a Gr\"obner basis in a lexicographic order induced
  by a vertex labeling. Such graphs are chordal and claw-free.  We give a reduced squarefree Gr\"obner basis for general
  $G$. It follows that all binomial edge ideals are radical ideals. Their minimal primes can be characterized by
  particular subsets of the vertices of $G$.  We provide sufficient conditions for Cohen--Macaulayness for closed and
  nonclosed graphs.

  Binomial edge ideals arise naturally in the study of conditional independence ideals.  Our results apply for the class
  of conditional independence ideals where a fixed binary variable is independent of a collection of other variables,
  given the remaining ones.  In this case the primary decomposition has a natural statistical interpretation.

\end{abstract}

\maketitle
\section*{Introduction}
Let $G$ be a simple graph on the vertex set $[n]=\{1,\ldots,n\}$, that is to say, $G$ has no loops and no multiple
edges. Furthermore let $K$ be a field and $S=K[x_1,\ldots,x_n,y_1,\ldots,y_n]$ be the polynomial ring in $2n$
variables. For $i<j$ we set $f_{ij}=x_iy_j-x_jy_i$. We define the {\em binomial edge} ideal $J_G\subset S$ of $G$ as the
ideal generated by the binomials $f_{ij}=x_iy_j-x_jy_i$ such that $i<j$ and $\{i,j\}$ is an edge of $G$. Note that if
$G$ has an isolated vertex $i$, and $G'$ is the restriction of $G$ to the vertex set $[n]\setminus \{i\}$, then
$J_G=J_{G'}$.

The class of binomial edge ideals is a natural generalization of the ideal of $2$-minors of a $2\times n$-matrix of
indeterminates. Indeed, the ideal of $2$-minors of a $2\times n$-matrix may be interpreted as the binomial edge ideal of
a complete graph on $[n]$. Related to binomial edge ideals are the ideals of adjacent minors considered by Ho\c{s}ten
and Sullivant \cite{HS}. In the case of a line graph our binomial edge ideal may be interpreted as an ideal of adjacent
minors. This particular class of binomial edge ideals has also been considered by Diaconis, Eisenbud and Sturmfels in
\cite{DES} where they compute the primary decomposition of this ideal.

Binomial edge ideals, as they are defined in this paper, also arise in the study of conditional independence statements
\cite{DSS}. They generalize a class which has been studied by Fink \cite{AF}.

Classically one studies edge ideals of a graph $G$ which are generated by the monomials $x_ix_j$ where $\{i,j\}$ is an
edge of $G$. The edge ideal of a graph has been introduced by Villarreal \cite{V} where he studied the Cohen--Macaulay
property of such ideals. The purpose of this paper is to study the algebraic properties of binomial edge ideals in terms
of properties of the underlying graph. In Section 1 we consider the Gr\"obner basis of $J_G$ with respect to the
lexicographic order induced by $x_1>x_2>\cdots >x_n>y_1>y_2>\cdots >y_n$. We show in Theorem~\ref{quadratic} that $J_G$
has a quadratic Gr\"obner basis if $G$ is closed with respect to the given labeling. Being closed can be characterized
by the associated acyclic directed graph $G^*$ with arrows $(i,j)$ whenever $\{i,j\}$ is an edge of $G$ and $i<j$. We
show in Proposition~\ref{characterization} that $G$ is closed if and only if for any two distinct vertices $i$ and $j$
of $G^*$, all shortest paths from $i$ to $j$ are directed. In Proposition~\ref{cmbinomial} we give a sufficient
condition for a closed graph to have a Cohen--Macaulay binomial edge ideal. In Theorem~\ref{mainresult} we compute
explicitly the reduced Gr\"obner basis of $J_G$ for any simple graph $G$.  This is one of the main results of this paper. As a
consequence we see that the initial ideal of $J_G$ is squarefree which in turn implies that $J_G$ is a reduced
ideal. Of course, Theorem~\ref{quadratic} is a simple consequence of  Theorem~\ref{mainresult}. But as the proof of Theorem~\ref{quadratic} is quite simple and as it leads to the concept of closed graphs, we decided to present Theorem~\ref{quadratic}  independent from Theorem~\ref{quadratic}.

Section 3 is devoted to the study of the minimal prime ideals of $J_G$. In Theorem~\ref{intersection} we write $J_G$ as a finite intersection of prime ideals which allows us to
compute the dimension of $S/J_G$. It turns out that if $S/J_G$ is Cohen--Macaulay, then $\dim S/J_G=|V(G)|+c$, where $c$
is the number of connected components of $G$. As a simple consequence of this, one sees that a circle of length $n$ is
unmixed or Cohen--Macaulay, if and only if $n=3$. As a last result of Section 3 we identify in Corollary\ref{minimalprime} the minimal prime ideals of $J_G$.  They are related to the cut-points of certain subgraphs of $G$.

In the last section we discuss applications to the study of
conditional independence ideals. For a class of conditional independence statements, suitable to model a notion of
robustness, the results in the prior sections show that the corresponding ideal is a radical ideal.  Furthermore, the
primary decomposition can be computed, which yields a classification and parametrization of the set of probability
distributions which satisfy these statements.

Terai informed the authors that M. Ohtani \cite{Oh} independently obtained similar results for this class of ideals.

\section{Edge ideals with quadratic Gr\"obner bases and closed graphs}

We first study  the question when $J_G$ has a quadratic Gr\"obner basis.

\begin{Theorem}
\label{quadratic}
Let $G$ be a simple graph on the vertex set $[n]$, and let $<$ be the lexicographic order on $S=K[x_1,\ldots,x_n,y_1,\ldots,y_n]$ induced by $x_1>x_2>\cdots >x_n>y_1>y_2>\cdots >y_n$. Then the  following conditions are equivalent:
\begin{enumerate}
\item[{\em (a)}] The generators $f_{ij}$ of $J_G$ form  a quadratic Gr\"obner basis;
\item[{\em (b)}] For all edges   $\{i,j\}$ and $\{k,l\}$ with $i<j$ and $k<l$ one has  $\{j,l\}\in E(G)$ if $i=k$, and $\{i,k\}\in E(G)$ if $j=l$.
\end{enumerate}
\end{Theorem}

\begin{proof}
(a) $\Rightarrow$ (b): Suppose (b) is violated, say, $\{i,j\}$ and $\{i,k\}$ are edges with $i<j<k$, but $\{j,k\}$ is not an edge. Then $S(f_{ik},f_{ij})=y_if_{jk}$ belongs to $J_G$, but none of the initial monomials of the  quadratic generators of $J_G$  divides $\ini_<(y_if_{jk})$.

(b) $\Rightarrow$ (a): We apply Buchberger's criterion and show that all $S$-pairs $S(f_{ij},f_{kl})$ reduce to  $0$. If $i\neq k$ and $j\neq l$, then $\ini_<(f_{ij})$  and $\ini_<(f_{kl})$ have no common factor. It is well known that in this case $S(f_{ij},f_{kl})$ reduces to  zero.
On the other hand, if $i=k$, we may assume that $l<j$. Then
\[
S(f_{ij},f_{il})=y_if_{lj}
\]
is the standard expression of $S(f_{ij},f_{il})$. Similarly, if $j=l$, we may assume that $i<k$. Then
\[
S(f_{ij},f_{kj})=x_jf_{ik}
\]
is the standard expression of $S(f_{ij},f_{kj})$. In both cases the $S$-pair reduces to  $0$.
\, \, \, \,
\end{proof}

Condition (b) of Theorem~\ref{quadratic} does not only depend on  the isomorphism type of the graph, but also on the labeling of its vertices. For example the  graph $G$ with edges $\{1,2\}$,  $\{2,3\}$, and  the graph $G'$ with edges $\{1,2\}$, $\{1, 3\}$ are isomorphic, but $G$ satisfies condition (b), while $G'$ does not.

In fact, condition (b) is a condition of the {\em associated directed graph} $G^*$ of $G$ which is defined as follows: the ordered pair $(i,j)$ is an arrow of $G^*$ if $\{i,j\}$ is an edge of $G$ with $i<j$. The directed graph $G^*$ is {\em acyclic}, that is, it has no directed cycles. Therefore we call $G^*$ also the associated acyclic directed graph of $G$.

An acyclic directed graph is also called an {\em acyclic digraph} or simply a {\em DAG}. Acyclic directed graphs constitute  an important class of directed graphs and play an important role in the modeling of information flows in networks. Any acyclic  directed graph arises in the same  way as we obtained $G^*$ from $G$. Indeed, one of the fundamental results on acyclic directed graphs $G$ is that they admit an {\em acyclic ordering} of its vertices, that is, the vertices of $G$ can be ordered $v_1,\ldots, v_r$ such that  for every arrow $(v_i,v_j)$ of $G$ we have $i<j$, see for example \cite[Proposition 1.4.3]{BG}. An acyclic directed graph usually has many different acyclic orderings. In \cite[Corollary 1.3]{St} Stanley expressed the number of possible acyclic orderings in terms of the chromatic polynomial of $G$.

We say that a graph $G$ on $[n]$ is {\em closed with respect to the given labeling of the vertices}, if $G$ satisfies condition (b) of  Theorem~\ref{quadratic}, and we say that a graph $G$ with vertex set $V(G)=\{v_1,\ldots,v_n\}$ is {\em closed}, if its vertices can be labeled by the integer $1,2,\ldots,n$ such that for this labeling $G$ is closed.

\begin{Proposition}
\label{hopefully} If $G$ is closed, then  $G$ is chordal and  has no induced  subgraph consisting of three different  edges  $e_1$,  $e_2$, $e_3$ with $e_1\sect e_2 \sect e_3\neq \emptyset$.
\end{Proposition}

\begin{proof} Suppose $G$ is not chordal, then $G$ contains a cycle $C$ of length $> 3$ with no chord. Let $i$ be the vertex of $C$ with $i<j$ for all $j\in V(C)$, and let $\{i,j\}$ and $\{i,k\}$ be the edges of $C$ containing $i$. Then $i<j$ and $i<k$, but $\{j,k\}\not\in E(G)$.

Since $G$ is closed, any induced subgraph is closed as well. Suppose there exists an induced subgraph $H$ with three different  edges  $e_1$,  $e_2$, $e_3$ such that three different  edges  $e_1$,  $e_2$, $e_3$ with $e_1\sect e_2 \sect e_3\neq \emptyset$. Then there exists $i$ such that  $e_1\sect e_2 \sect e_3=\{i\}$. Say, $e_1=\{i,j\}$, $e_2=\{i,k\}$ and  $e_3=\{i,l\}$. Then $i\neq \min\{i,j,k,l\}$, otherwise $H$ is not closed. If $j<i$, then $k>i$ and $l>i$,  since $H$ is closed. But then $\{k,j\}$ must be an edge of $H$, a contradiction.
\, \, \, \, \, \, \, \, \, \,
\, \, \, \, \, \, \, \, \, \,
\, \, \, \, \, \, \, \, \, \,
\, \, \, \, \, \, \, \, \, \,
\, \, \, \, \, \, \, \, \, \,
\end{proof}

A graph with three different  edges  $e_1$,  $e_2$, $e_3$ such that  $e_1\sect e_2 \sect e_3\neq \emptyset$ is called a {\em claw}. Hence Proposition~\ref{hopefully} says that a closed graph is a claw-free chordal graph.

\begin{Corollary}
\label{bipartite}
A  bipartite graph is closed if and only if it is a line.
\end{Corollary}

\begin{proof} A bipartite graph has no odd cycles. Since a closed graph is chordal, and since a chordal graph has an odd cycle, unless it is a tree, a closed bipartite graph must be a tree. If the tree is not a line, then there exists an induced subgraph which is a claw. Thus a closed bipartite graph must be a line.

Conversely, if $G$ is a line of length $l$, then $G$ is closed for the   labeling of the vertices such that  $\{1,2\},\{2,3\},\ldots,\{l,l+1\}$ are the  edges of $G$.
\, \, \, \, \, \, \, \, \, \,
\, \, \, \, \, \, \, \, \, \,
\end{proof}

The conditions for being a closed graph formulated in Proposition~\ref{hopefully} are only sufficient. For example the graph with edges $\{a,b\}$, $\{b,c\}$, $\{a,c\}$, $\{a,x\}$,$\{b,y\}$ and $\{c,z\}$ is chordal without a claw, but is not closed.

In the following we give a characterization of graphs which are closed with respect to a given labeling. Let $G$ be a graph, and let $v$ and $w$ be vertices of $G$. A {\em path} $\pi$ from $v$ to $w$ is a sequence of vertices $v=v_0,v_1,\ldots,v_l=w$ such that each $\{v_i,v_{i+1}\}$ is an edge of the underlying graph. If $G$ is directed, then the path $\pi$ is called {\em directed}, if either $(v_i,v_{i+1})$ is an arrow for all $i$, or $(v_{i+1},v_i)$ is an arrow for all $i$.

\begin{Proposition}
\label{characterization}
A graph $G$  on $[n]$ is closed with respect to the given labeling, if and only if  for any two vertices $i\neq j$ of associated directed graph $G^*$, all paths of shortest length  from $i$ to $j$ are directed.
\end{Proposition}

\begin{proof}
Suppose all shortest paths from $i$ to $j$ in $G^*$ are directed. Let $(i,j)$ and $(i,k)$ be two arrow with $j<k$.  Then $\{j,i\},\{i,k\}$ is a path from $j$ to $k$ which is not directed. So it cannot be the shortest path. Hence there exists the arrow $(j,k)$. Similarly it follows that if $(i,k)$  and $(j,k)$  are arrows of $G^*$ with $i<j$, then there must exist the arrow $(i,j)$ in $G^*$. This shows that $G^*$ is closed.

Conversely, assume that $G$ is closed.  Then there exists a labeling such that $G^*$ is closed. Let $i$ and $j$ be two distinct vertices and let $P$ be path of shortest length from $i$ to $j$. Suppose $P$ is not directed. Then there there exists a subpath $r,s,t$ of $P$  such that  $(r,s)$, $(t,s)$,  or  $(s,r),s(s,t)$ in $G^*$. In both cases we may assume that $r<t$. Then, since $G^*$ is closed, it follows that $(r,t)$  is an arrow in $G^*$. Replacing the subpath $r,s,t$   by $r,t$,  we obtain a shorter path from $i$ to $j$, a contradiction.
\, \, \, \, \, \, \, \, \, \,
\, \, \, \, \, \, \, \, \, \,
\, \,
\end{proof}

In Proposition~\ref{characterization} it is important to  require that {\em all} paths of shortest length from $i$ to $j$ are directed in order to conclude that $G^*$ is closed. Indeed, consider the graph $G$ with edges  $\{1,2\}$, $\{2,3\}$, $\{3,4\}$ and $\{1,4\}$. Then the path $2,3,4$ is directed, while $2,1,4$ is not directed. But both paths are shortest paths between $2$ and $4$.

\begin{Proposition}
\label{closure}
Let $G$ be a simple graph on $[n]$. Then there exists a unique minimal (with respect to inclusion of edges)  graph $\bar{G}$  on $[n]$ whose associated acyclic graph is closed  with respect to the given labeling and  such that $G$ is a subgraph of $\bar{G}$.
\end{Proposition}

\begin{proof}
 Consider the set $\mathcal{C}$  of  graphs on $[n]$ containing $G$ and whose associated acyclic graph is closed. This set is  not empty, because the complete graph on $[n]$ belongs to this set. Since the intersection of any two  graphs in $\mathcal{C}$  belongs again to $\mathcal{C}$, the assertion follows, as desired.
\, \, \, \, \, \, \, \, \, \,
\, \, \, \, \, \, \, \, \, \,
\, \, \, \, \, \, \, \, \, \,
\, \, \, \, \, \, \, \, \, \,
\, \, \, \, \, \, \, \, \, \,
\, \, \, \,
\end{proof}

The unique  minimal closed graph $\bar{G}$ containing $G$ is called the {\em closure} of $G$.

\medskip
One  basic question is which of the binomial edge ideals  are Cohen--Macaulay. For a graph $G$,  this is the case if and only the binomial edge of each component is Cohen--Macaulay. Thus it is enough to consider connected graphs.   A partial answer on  the Cohen--Macaulayness of binomial edge ideals  is given in

\begin{Proposition}
\label{cmbinomial}
Let $G$ be a  connected graph on $[n]$ which  is closed with respect to the given labeling. Suppose further that  $G$ satisfies the condition that whenever  $\{i,j+1\}$ with $i<j$ and $\{j,k+1\}$ with $j<k$ are edges of $G$,  then $\{i,k+1\}$ is an edge of $G$. Then $S/J_G$ is Cohen--Macaulay.
\end{Proposition}

\begin{proof}
We will show that $S/\ini_<(J_G)$ is Cohen--Macaulay. This will then imply that $S/J_G$ is Cohen--Macaulay as well.

Since the associated acyclic directed graph is closed, it follows from Theorem~\ref{quadratic} that $\ini_<(J_G)$ is generated by the monomials $x_iy_j$ with $\{i,j\}\in E(G)$ and $i<j$. Applying the automorphism $\varphi \: S\to S$ which maps each $x_i$ to $x_i$,  and $y_j$ to $y_{j-1}$ for $j>1$ and $y_1$ to $y_n$,  $\ini_<(J_G)$ is mapped to the ideal  generated by all monomials $x_iy_j$ with $\{i,j+1\}\in E(G)$. This ideal  has all its generators in $S'=K[x_1,\ldots,x_{n-1},y_1,\ldots,y_{n-1}]$. Let $I\subset S'$ be the ideal generated by these monomials. Then $S/\ini_<(J_G)$ is Cohen--Macaulay if and only if $S'/I$ is Cohen--Macaulay. Note that $I$ is the edge ideal of the bipartite graph $\Gamma$ on the vertex set $\{x_1,\ldots,x_{n-1},y_1,\ldots,y_{n-1}\}$, and  with $\{x_i,y_j\}\in E(\Gamma)$ if and only if $\{i,j+1\}\in E(G)$. In \cite{HH} the Cohen--Macaulay bipartite graphs are characterized as follows: Suppose the edges of the bipartite graph can be labeled such that
\begin{enumerate}
\item[(i)] $\{x_i,y_i\}$ are edges for $i=1,\ldots,n$;

\item[(ii)] if $\{x_i,y_j\}$ is an edge, then $i\leq j$;

\item[(iii)] if $\{x_i,y_j\}$ and $\{x_j,y_k\}$ are edges, then
$\{x_i,y_k\}$ is an edge.
\end{enumerate}
Then the corresponding edge ideal is Cohen--Macaulay.

We are going to verify these conditions for our edge ideal. Condition (ii) is trivially satisfied, and condition (iii) is a consequence of our assumption that  whenever  $\{i,j+1\}$ with $i<j$ and $\{j,k+1\}$ with $j<k$ are edges of $G$,  then $\{i,k+1\}$ is an edge of $G$.

For condition (i) we have to show that $\{i,i+1\}\in E(G)$ for all $i$. But this follows from Proposition~\ref{characterization} which says that all shortest paths from $i$ to $i+1$ are oriented paths. If $i,i+1$ would not be a path, then a shortest path from $i$ to $i+1$ could not be oriented. Thus $i,i+1$ is a path in $G$, and hence $\{i,i+1\}\in E(G)$.
\, \, \, \, \, \, \, \, \, \,
\, \, \, \, \, \, \, \, \, \,
\end{proof}

\begin{Examples}
\label{cmexamples}
{\em (a) Any complete graph satisfies the conditions  of Proposition~\ref{cmbinomial}, so that $S/J_G$ is Cohen--Macaulay. But of course this is well known because in this case $J_G$ is the ideal of $2$-minors of a generic $2\times n$-matrix.

(b) Any line graph with the natural order of the vertices satisfies conditions  of Proposition~\ref{cmbinomial}. Actually $J_G$ is a complete intersection in this case.

(c) There are many more graphs satisfying the conditions of Proposition~\ref{cmbinomial}. For example the graph with edges $\{1,2\}$, $\{2,3\}$ $\{1,3\}$ and $\{3,4\}$.

(d) Not all closed graphs satisfy the conditions of Proposition~\ref{cmbinomial}. Such an example is the graph with edges $\{1,2\}$, $\{1,3\}$,  $\{2,3\}$, $\{1,4\}$ and $\{3,4\}$. For this graph we have that  $\ini_<(J_G)$ and $J_G$ are  not Cohen--Macaulay.

(e) A graph $G$ need not be closed for $S/J_G$ being Cohen--Macaulay. The  graph given after Corollary \ref{bipartite}  is such an example. }
\end{Examples}

\section{The reduced Gr\"obner basis of a binomial edge ideal}
We now come to the main result of this paper.
For this we need to introduce the following concept:
let $G$ be a simple graph on $[n]$,
and let $i$ and $j$ be two vertices of $G$ with $i<j$.
A path $i=i_0,i_1,\ldots,i_r=j$ from $i$ to $j$
is called {\em admissible}, if
\begin{enumerate}
\item[(i)] $i_k\neq i_\ell$  for $k\neq \ell$;
\item[(ii)] for each $k=1,\ldots,r-1$ one has either $i_k<i$ or $i_k>j$;
\item[(iii)] for any proper subset $\{j_1,\ldots,j_s\}$
of $\{i_1,\ldots,i_{r-1}\}$, the sequence $i,j_1,\ldots,j_s,j$
is not a path.
\end{enumerate}
Given an admissible path
\[
\pi: i=i_0,i_1,\ldots,i_r=j
\]
from $i$ to $j$, where $i < j$, we associate the monomial
\[
u_{\pi}=(\prod_{i_k>j}x_{i_k}) (\prod_{i_\ell<i}y_{i_\ell}).
\]

\begin{Theorem}
\label{mainresult}
Let $G$ be a simple graph on $[n]$. Let $<$ be the monomial order introduced in Theorem~\ref{quadratic}. Then the set of binomials
\[
{\mathcal G}
= \Union_{i<j} \,
\{\,u_{\pi}f_{ij}\,:\;\text{$\pi$ is an admissible path from $i$ to $j$}\,\}
\]
is a reduced Gr\"obner basis of $J_G$.
\end{Theorem}

\begin{proof}
We organize this proof as follows:
In First Step, we prove that ${\mathcal G} \subset J_G$.
Then, since ${\mathcal G}$ is a system of generators,
in Second Step, we show that ${\mathcal G}$ is a Gr\"obner
basis of $J_G$ by using Buchberger's criterion.  Finally,
in Third Step, it is proved that ${\mathcal G}$ is reduced.

\medskip
\noindent
{\em First Step.}
We show that, for each admissible path $\pi$ from $i$ to $j$,
where $i < j$,
the binomial $u_\pi f_{ij}$ belongs $J_G$.
Let
$\pi: i=i_0,i_1,\ldots,i_{r-1},i_r=j$ be an admissible path in $G$.
We proceed with induction on $r$.
Clearly the assertion is true if $r = 1$.
Let $r > 1$ and
$A = \{ i_k : i_k < i \}$
and
$B = \{ i_\ell : i_\ell > j \}$.
One has either $A \neq \emptyset$ or $B \neq \emptyset$.
If $A \neq \emptyset$, then we set $i_{k_0} = \max A$.
If $B \neq \emptyset$, then we set $i_{\ell_0} = \min B$.

Suppose $A \neq \emptyset$.
It then follows that
each of the paths
$\pi_1 : i_{k_0}, i_{k_0-1}, \ldots, i_1, i_0=i$
and
$\pi_2 : i_{k_0}, i_{k_0+1}, \ldots, i_{r-1}, i_r = j$
in $G$ is admissible.
Now,
the induction hypothesis guarantees that
each of $u_{\pi_1}f_{i_{k_0},i}$ and
$u_{\pi_2}f_{i_{k_0},j}$ belongs to $J_G$.
A routine computation says that the $S$-polynomial
$S(u_{\pi_1}f_{i_{k_0},i}, u_{\pi_2}f_{i_{k_0},j})$
is equal to $u_\pi f_{ij}$.
Hence $u_\pi f_{ij} \in J_G$, as desired.

When $B \neq \emptyset$, the same argument as
in the case $A \neq \emptyset$ is valid.

\medskip
\noindent
{\em Second Step.}
It will be proven that the set of those
binomials $u_\pi f_{ij}$, where $\pi$ is an admissible path from $i$ to $j$,
forms a Gr\"obner basis of $J_G$.
In order to show this we apply Buchberger's criterion, that is,
we show that all $S$-pairs $S(u_\pi f_{ij}, u_\sigma f_{k\ell})$,
where $i < j$ and $k < \ell$,
reduce to zero. For this we will consider different cases.

In the case that $i=k$ and $j=\ell$, one has $S(u_\pi f_{ij}, u_\sigma f_{k\ell})=0$.

In the case that $\{i,j\}\sect\{k,\ell\}=\emptyset$, or $i=\ell$, or $k=j$, the initial monomials $\ini_<(f_{ij})$ and $\ini_<(f_{k\ell})$
form a regular sequence. Hence the $S$-pair $S(u_\pi f_{ij}, u_\sigma f_{k\ell})$ reduce to zero, because of the following more general fact: let $f,g\in S$ such that $\ini_<(f)$ and $\ini_<(g)$ form a regular sequence and let $u$ and $v$ be any monomials. Then $S(uf,vg)$ reduces to zero.

It remains to consider the cases that either $i=k$ and $j\neq \ell$ or $i\neq k$ and $j=\ell$.
Suppose we are in the first case.
(The second case can be proved similarly.)
We must show that $S(u_\pi f_{ij}, u_\sigma f_{i\ell})$ reduces to zero. We may assume that $j<\ell$, and must find a standard expression for $S(u_\pi f_{ij}, u_\sigma f_{i\ell})$
whose remainder is equal to zero.

Let $\pi: i=i_0,i_1,\ldots,i_r=j$ and $\sigma\: i=i_0',i_1',\ldots,i_s'=\ell$. Then there exist unique indices $a$ and $b$ such that
\[
i_a=i_b'\quad \text{and}\quad \{i_{a+1},\ldots,i_r\}\sect\{i_{b+1}',\ldots,i_s'\}=\emptyset.
\]
Consider the path
\[
\tau: j=i_r,i_{r-1},\ldots,i_{a+1},i_a=i'_b,i'_{b+1},\ldots,i'_{s-1}, i'_s=\ell
\]
from $j$ to $\ell$. To simplify the notation we write this path as
\[
\tau\: j=j_0,j_1,\ldots,j_t=\ell.
\]
Let
\[
j_{t(1)} = \min\{ \, j_c \, : \, j_c > j, \, c = 1, \ldots, t \, \},
\]
and
\[
j_{t(2)} = \min\{ \, j_c \, : \, j_c > j, \, c = t(1) + 1, \ldots, t \, \}.
\]
Continuing these procedures yield the integers
\[
0 = t(0) < t(1) < \cdots < t(q-1) < t(q) = t.
\]
It then follows that
\[
j =j _{t(0)} <  j_{t(1)} < \cdots < j_{t(q)-1} < j_{t(q)} = \ell
\]
and, for each $1 \leq c \leq t$, the path
\[
\tau_c : j_{t(c-1)}, j_{t(c-1)+1}, \ldots, j_{t(c)-1}, j_{t(c)}
\]
is admissible.

The highlight of the proof is to show that
\[
S(u_\pi f_{ij}, u_\sigma f_{i\ell})
= \sum_{c=1}^{q} v_{\tau_c} u_{\tau_c}f_{j_{t(c-1)} j_{t(c)}}
\]
is a standard expression of $S(u_\pi f_{ij}, u_\sigma f_{i\ell})$
whose remainder is equal to $0$,
where each $v_{\tau_c}$ is the monomial defined as follows:
Let $w = y_i \lcm(u_\pi,u_\sigma)$.  Thus $S(u_\pi f_{ij}, u_\sigma f_{i\ell}) = - w f_{j\ell}$.  Then
\begin{enumerate}
\item[(i)] if $c = 1$, then
\[
v_{\tau_1}=\frac{x_{\ell}w}{u_{\tau_1}x_{j_{t(1)}}};
\]
\item[(ii)] if $1 < c < q$, then
\[
v_{\tau_c}=\frac{x_{j}x_{\ell}w}{u_{\tau_c} x_{j_{t(c-1)}}x_{j_{t(c)}}};
\]
\item[(iii)] if $c = q$, then
\[
v_{\tau_q}=\frac{x_{j}w}{u_{\tau_q}x_{j_{t(q-1)}}}.
\]
\end{enumerate}
Our work is to show that
\[
w f_{j\ell}
= \frac{w x_{\ell}}{x_{j_{t(1)}}} f_{j j_{t(1)}} +
\sum_{c=2}^{q-1} \frac{w x_{j}x_{\ell}}{x_{j_{t(c-1)}}x_{j_{t(c)}}}
f_{j_{t(c-1)} j_{t(c)}}
+\frac{w x_{j}}{x_{j_{t(q-1)}}} f_{j_{t(q-1)} \ell}
\]
is a standard expression of $w f_{j\ell}$ with remainder $0$.
In other words, we must prove that
\begin{eqnarray*}
(\sharp)
\, \, \, \, \, \, \, \, \, \,
\, \, \, \, \, \, \, \, \, \,
w (x_jy_\ell - x_\ell y_j)
& = &\frac{w x_{\ell}}{x_{j_{t(1)}}} (x_{j}y_{j_{t(1)}} - x_{j_{t(1)}}y_{j}) \\
&   & + \,\,\, \sum_{c=2}^{q-1} \frac{w x_{j}x_{\ell}}{x_{j_{t(c-1)}}x_{j_{t(c)}}}
(x_{j_{t(c-1)}}y_{j_{t(c)}} - x_{j_{t(c)}}y_{j_{t(c-1)}}) \\
&   & + \,\,\, \frac{w x_{j}}{x_{j_{t(q-1)}}} (x_{j_{t(q-1)}}y_{\ell} - x_{\ell}y_{j_{t(q-1)}})
\end{eqnarray*}
is a standard expression of $w (x_jy_\ell - x_\ell y_j)$ with remainder $0$.

Since
\begin{eqnarray*}
w x_jy_\ell \, = \, \frac{w x_{j}}{x_{j_{t(q-1)}}}x_{j_{t(q-1)}}y_{\ell}
& > & \frac{w x_{j}x_{\ell}}{x_{j_{t(q-2)}}x_{j_{t(q-1)}}}
x_{j_{t(q-2)}}y_{j_{t(q-1)}} \\
& > & \cdots \,
\, > \, \frac{w x_{j}x_{\ell}}{x_{j_{t(1)}}x_{j_{t(2)}}} x_{j_{t(1)}}y_{j_{t(2)}}
\, > \, \frac{w x_{\ell}}{x_{j_{t(1)}}} x_{j}y_{j_{t(1)}},
\end{eqnarray*}
it follows that, if the equality $(\sharp)$ holds, then
$(\sharp)$ turns out to be
a standard expression of $w (x_j y_\ell - x_\ell y_j)$ with remainder $0$.
If we rewrite $(\sharp)$ as
\begin{eqnarray*}
\, \, \, \, \, \, \, \, \, \,
\, \, \, \, \, \, \, \, \, \,
w (x_jy_\ell - x_\ell y_j)
& = & w ( x_j x_\ell \frac{y_{j_{t(1)}}}{x_{j_{t(1)}}} - x_\ell y_j )
\\
&   & + \,\,\, w x_j x_\ell \sum_{c=2}^{q-1}
(
\frac{y_{j_{t(c)}}}{x_{j_{t(c)}}}
-
\frac{y_{j_{t(c-1)}}}{x_{j_{t(c-1)}}}
) \\
&   & + \,\,\,
w (
x_j y_\ell -
x_j x_\ell \frac{y_{j_{t(q-1)}}}
{x_{j_{t(q-1)}}}
) ,
\end{eqnarray*}
then clearly the equality holds.

\medskip
\noindent
{\em Third Step.}
Finally, we show that the Gr\"obner basis
${\mathcal G}$ is reduced.
Let $u_\pi f_{ij}$ and $u_\sigma f_{k\ell}$,
where $i < j$ and $k < \ell$, belong to
${\mathcal G}$ with $u_\pi f_{ij} \neq u_\sigma f_{k\ell}$.
Let $\pi : i = i_0, i_1, \ldots, i_r = j$
and $\sigma : k = k_0, k_1, \ldots, k_s = \ell$.
Suppose that $u_\pi x_iy_j$ divides either
$u_\sigma x_ky_\ell$ or $u_\sigma x_\ell y_k$.
Then $\{ i_0, i_1, \ldots, i_r \}$ is a proper subset
of $\{ k_0, k_1, \ldots, k_s \}$.

Let $i = k$ and $j = \ell$.
Then $\{ i_1, \ldots, i_{r-1} \}$ is a proper subset
of $\{ k_0, k_1, \ldots, k_s \}$
and $k, i_1, \ldots, i_{r-1}, \ell$ is an admissible path.
This contradicts the fact that $\sigma$ is an admissible path.

Let $i = k$ and $j \neq \ell$.
Then $y_j$ divide $u_\sigma$.
Hence $j < k$.  This contradicts $i < j$.

Let $\{ i, j \} \cap \{ k, \ell \} = \emptyset$.
Then $x_iy_j$ divide $u_\sigma$.
Hence $i > \ell$ and $j < k$.
This contradicts $i < j$.
\, \, \, \, \, \, \, \, \, \,
\, \, \, \, \, \, \, \, \, \,
\, \, \, \, \, \, \, \, \, \,
\, \, \, \, \, \, \, \, \, \,
\, \, \, \, \, \, \, \, \, \,
\, \, \, \, \, \,
\end{proof}

\begin{Corollary}
\label{radical}
$J_G$ is a radical ideal.
\end{Corollary}

\begin{proof}
The assertion follows from Theorem~\ref{mainresult} and the following general fact: let $I\subset S$ be a graded ideal with the property that $\ini_<(I)$ is squarefree for some monomial order $<$. Then $I$ is a radical ideal. Indeed, there exists an ideal $\tilde{I}\subset S[t]$ in the polynomial ring $S[t]$   such that $t$ is a nonzerodivisor on $S[t]/\tilde{I}$ with $(S[t]/\tilde{I})/(tS[t]/\tilde{I})\iso S/\ini_<(I)$ and such  that $\tilde{I}S[t,t^{-1}]=IS[t,t^{-1}]$,  and there are positive degrees  on the variables of $K[x_1,\ldots,x_n, t]$ such that $\tilde{I}$ is a graded ideal with respect to this grading. Thus we may apply the graded version of Lemma 4.4.9 in \cite{BH}  in order to conclude that $\tilde{I}$ is a radical ideal. From the equality $\tilde{I}S[t,t^{-1}]=IS[t,t^{-1}]$, it follows that $I$ is a radical ideal as well.
\, \, \, \, \, \, \, \,
\end{proof}

As a consequence of Theorem~\ref{mainresult}  we see that  all admissible paths of a graph $G$ can be determined by computing the reduced Gr\"obner basis of $J_G$.

On the other hand, it is not the case that for each edge $\{i,j\}$ in the closure of $G$ there exists an admissible path from $i$ to $j$. For example, for the graph $G$  with edges $\{2,3\}$, $\{1,3\}$ and $\{1,4\}$, the edge $\{2,4\}$ belongs to the closure of $G$, but the only path $2,3,1,4$ from $2$ to $4$ is not admissible.
Thus the reduced Gr\"obner basis of $J_G$ does not give the closure of $G$.

\section{The minimal prime ideals of a binomial edge ideal}

Let $G$ be a simple graph on $[n]$. For each subset $S\subset [n]$ we define a prime ideal $P_S$. Let $T=[n]\setminus S$, and let $G_1,\ldots, G_{c(S)}$ be the connected component of $G_T$. Here $G_T$ is the {\em restriction} of $G$ to $T$ whose edges are exactly those edges $\{i,j\}$ of $G$ for which $i,j\in T$. For each $G_i$ we denote by $\tilde{G}_i$ the complete graph on the vertex set $V(G_i)$. We set
\[
P_S(G)=(\Union_{i\in S}\{x_i,y_i\}, J_{\tilde{G_1}},\ldots, J_{\tilde{G}_{c(S)}}).
\]
Obviously, $P_S(G)$ is a prime ideal. In fact, each $J_{\tilde{G}_i}$ is the ideal of $2$-minors of a  generic $2\times n_j$-matrix with $n_j=|V(G_j)|$. Since all the prime ideals $J_{\tilde{G}_j}$, as well as the prime ideal $(\Union_{i\in S}\{x_i,y_i\})$ are prime ideals in  pairwise different sets of variables,  $P_S(G)$ is a prime ideal, too.

\begin{Lemma}
\label{height}
With the notation introduced we have $\height P_S(G)=|S|+(n-c(S))$.
\end{Lemma}

\begin{proof}
We have
\begin{eqnarray*}
\height P_S(G)&=&\height (\Union_{i\in S}\{x_i,y_i\})+\sum_{j=1}^{c(S)}\height  J_{\tilde{G_j}}=2|S|+\sum_{j=1}^{c(S)}(n_j-1)\\
&=& |S|+(|S|+\sum_{j=1}^{c(S)}n_j) -c(S)=|S|+(n-c(S)),
\end{eqnarray*}
as required.
\, \, \, \, \, \, \, \, \, \,
\, \, \, \, \, \, \, \, \, \,
\, \, \, \, \, \, \, \, \, \,
\, \, \, \, \, \, \, \, \, \,
\, \, \, \, \, \, \, \, \, \,
\end{proof}

In \cite{ES} Eisenbud and Sturmfels showed that all associated prime ideals of a binomial ideal are binomial ideals. In our particular case we have

\begin{Theorem}
\label{intersection}
Let $G$ be a simple graph on the vertex set $[n]$. Then $J_G=\Sect_{S\subset [n]}P_S(G)$.
\end{Theorem}

\begin{proof}
It is obvious that each of the prime ideals $P_S(G)$ contains $J_G$.
We will show by induction on $n$ that each minimal prime ideal containing $J_G$ is of the form $P_S(G)$ for some $S\subset [n]$. Since by Corollary~\ref{radical},  $J_G$ is a radical ideal, and since a radical ideal is the intersection of its minimal prime ideals, the assertion of the theorem will follow.

We may assume that $G$ is connected. Because if $G_1,\ldots, G_r$ are the connected components of $G$, then each minimal prime ideal $P$ of $J_G$ is of the form $P_1+\ldots +P_r$ where each $P_i$ is a minimal prime ideal of $J_{G_i}$. Thus if each $P_i$ has the expected form, then so does $P$.
So now let $G$ be connected and let  $P$ be a  minimal prime ideal of $J_G$.  Let $T$ be the maximal subset of $\{x_1,\ldots,x_n\}$ with the property that  $T\subset P$ and that  $x_i\in T$ implies $y_i\not\in P$. We will show that $T=\emptyset$. This will then imply that if $x_i\in P$, then $y_i\in P$, as well.

We first observe that $T\neq \{x_1,\ldots,x_n\}$. Because otherwise we would have  $J_G\subset J_{\tilde{G}}\subsetneq (x_1,\ldots,x_n)\subset P$, and $P$ would not be a minimal prime ideal of $J_G$.

Suppose that $T\neq \emptyset$. Since $T\neq\{x_1,\ldots,x_n\}$, and since $G$ is connected there exists $\{i,j\}\in E(G)$ such that $x_i\in T$ but $x_j\not\in T$. Since $x_iy_j-x_jy_i\in J_G\subset P$, and since $x_i\in P$ it follows that $x_jy_i\in P$.  Hence since $P$ is a prime ideal, we have $x_j\in P$ or $y_i\in P$. By the definition of $T$ the second case cannot happen, and so $x_j\in P$. Since $x_j\not\in T$,  it follows that $y_j\in P$.

Let $G'$ be the restriction of $G$ to the vertex set to $[n]\setminus \{j\}$. Then
\[
(J_{G'},x_j,y_j)=(J_G, x_j, y_j)\subset P.
\]
Thus $\bar{P}=P/(x_j,y_j)$ is a minimal prime ideal of $J_{G'}$ with $x_i\in \bar{P}$ but $y_i\not \in \bar{P}$ for all $x_i\in T\subset \bar{P}$.  By induction hypothesis, $\bar{P}$ is of the form $P_S(G')$ for some subset $S\subset [n]\setminus\{j\}$. This contradicts the fact that $T\neq \emptyset$.

By what we have shown it follows that there exists a  subset $S\subset[n]$ such that $P=(\Union_{i\in S}\{x_i,y_i\}, \bar{P})$ where $\bar{P}$ is a prime ideal containing no variables. Let $G'$ be the graph $G_{[n]\setminus S}$. Then reduction modulo the ideal $(\Union_{i\in S}\{x_i,y_i\})$ shows that $\bar{P}$ is a monomial prime ideal $J_{G'}$ which contains no variables. Let $G_1, \ldots, G_c$ be the connected components of $G'$. We will show that $\bar{P}=(J_{\tilde{G_1}},\ldots, J_{\tilde{G}_{c}})$. This then implies that $P=(\Union_{i\in S}\{x_i,y_i\}, J_{\tilde{G_1}},\ldots, J_{\tilde{G}_{c}})$, as desired.

To simplify notation we may as well assume that $P$ contains no variables and have to show that $P=(J_{\tilde{G_1}},\ldots, J_{\tilde{G}_{c}})$, where $G_1,\ldots,G_c$ are the connected components of $G$. In order to prove this we claim that if $i,j$ with $i<j$ are two edges of $G_k$ for some $k$, then $f_{ij}\in P$. From this it will then follow that $(J_{\tilde{G_1}},\ldots, J_{\tilde{G}_{c}})\subset P$. Since $(J_{\tilde{G_1}},\ldots, J_{\tilde{G}_{c}})$ is a prime ideal containing $J_G$, and $P$ is a minimal prime ideal containing $J_G$,  we conclude that  $P=(J_{\tilde{G_1}},\ldots, J_{\tilde{G}_{c}})$.

Let $i=i_0,i_1,\ldots,i_r=j$ a path in $G_k$ from $i$ to $j$. We proceed by induction on $r$ to show that $f_{ij}\in P$. The assertion is trivial for $r=1$. Suppose now that $r>1$. Our induction hypothesis says that $f_{i_1 j}\in P$. On the other hand, one has $x_{i_1}f_{ij}=x_jf_{ii_1}+x_if_{i_1j}$.  Thus $x_{i_1}f_{ij}\in P$. Since $P$ is a prime ideal
and since $x_{i_1}\not \in P$, we see that $f_{ij}\in P$.
\, \, \, \, \, \,
\end{proof}

Lemma~\ref{height} and Theorem~\ref{intersection} yield the following

\begin{Corollary}
\label{dimension}
Let $G$ be a simple graph on $[n]$. Then $$\dim S/J_G=\max\{(n-|S|)+c(S)\;\:\; S\subset [n]\}.$$
In particular,  $\dim S/J_G\geq n+c$, where $c$ is the number of connected components of $G$.
\end{Corollary}

In general, this inequality is strict. For example, for our claw $G$ with edges $\{1,2\}$, $\{1,3\}$ and $\{1,4\}$ we have $\dim S/J_G=6$.

\begin{Corollary}
\label{beingcm}
Let $G$ be a simple graph on $[n]$ with $c$ connected components. If $S/J_G$ is Cohen--Macaulay, then $\dim S/J_G=n+c$.
\end{Corollary}

\begin{proof}
Since $P_\emptyset(G)$ does not contain any monomials, it follows that $P_S(G)\nsubseteq P_\emptyset(G)$ for any nonempty subset $S\subset [n]$. Thus Theorem~\ref{intersection} implies that $P_\emptyset(G)$ is a minimal prime ideal of $J_G$. Since $\dim S/P_\emptyset(G)=n+c$
and since $S/J_G$ is equidimensional, the assertion follows.
\, \, \, \, \, \, \, \, \, \,
\, \, \, \, \, \, \, \, \, \,
\, \, \, \, \, \, \, \, \, \,
\, \, \, \, \, \, \, \, \, \,
\, \, \, \, \, \, \, \, \, \,
\, \, \, \, \,
\end{proof}

\begin{Example}
{\em
Consider the line graph $G$ with $n$ vertices. Then, as observed in Example \ref{cmexamples}, $S/J_G$ is Cohen--Macaulay. It follows from Corollary~\ref{beingcm} that $\dim S/P=n+1$ for all minimal prime ideals of $J_G$. Let $S$ be any subset of $[n]$. Then Theorem~\ref{intersection} and Corollary~\ref{dimension} imply that  the  minimal prime ideals of $J_G$  are exactly those prime ideals $P_S(G)$ for which $c(S)=|S|+1$. Let $S\subset [n]$.  Then there exists integers $1\leq a_1\leq b_1<a_2\leq b_2<a_3\leq b_3<\cdots <a_r\leq b_r\leq n$ such that
\[
S=\Union_{i=1}^r[a_i,b_i]\quad \text{where for each $i$,}\quad [a_i,b_i]=\{j\in\ZZ\: a_i\leq j\leq b_i \}.
\]
We see that $|S|=\sum_{i=1}^r(b_i-a_i+1)=\sum_{i=1}^r(a_i-b_i)+r$, and that
\[
c(S)= \left\{ \begin{array}{ll}
       r-1, & \;\text{if $a_1=1$ and $b_r=n$}, \\
       r, & \;\text{if $a_1\neq 1$ and $b_r=n$, or $a_1=1$ and $b_r\neq n$,}\\
       r+1, &\; \text{if $a_1\neq 1$ and $b_r\neq n$.}
        \end{array} \right.
\]
Thus $c(S)=|S|+1$ if and only if $a_1\neq 1$, $b_r\neq n$ and $a_i=b_i$ for all $i$. In other words, the minimal prime ideals of $G$ are those $P_S(G)$ for which $S$ is a subset of $[n]$ of the form $\{a_1,a_2,\ldots,a_r\}$ with $1<a_1<a_2<\ldots <a_r<n$. This is exactly the result of Diaconis, Eisenbud and  Sturmfels \cite[Theorem 4.3]{DES}.
}
\end{Example}

The question of when $J_G$ is a prime ideal is easy to answer.

\begin{Proposition}
\label{beingadomain}
Let $G$ be a simple graph on $[n]$. Then $J_G$ is a prime ideal if and only if each connected component of $G$ is a complete graph.
\end{Proposition}

\begin{proof}
Let $G_1,\ldots,G_r$ be the connected components of $G$, and suppose that $J_G$ is a prime ideal. Since $P_\emptyset(G)=(J_{\tilde{G_1}},\ldots, J_{\tilde{G_r}})$ is a minimal prime ideal of $J_G$ and $J_G$ is a prime ideal, it follows that $J_G=(J_{\tilde{G_1}},\ldots, J_{\tilde{G_r}})$. On the other hand, $J_G=(J_{G_1},\ldots, J_{G_r})$. Thus the desired conclusion is a consequence of the following observation. Suppose that  $G$ and $G'$ are graphs on $[n]$ with $V(G)\subset V(G')$. Then $E(G)=E(G')$, if and only $J_G=J_{G'}$.
\end{proof}

\begin{Corollary}
\label{cycle}
Let $G$ be a cycle of length $n$. Then the following conditions are equivalent:
\begin{enumerate}
\item[{\em (a)}] $n=3$.
\item[{\em (b)}] $J_G$ is a prime ideal.
\item[{\em (c)}] $J_G$ is unmixed.
\item[{\em (d)}] $S/J_G$ is Cohen--Macaulay.
\end{enumerate}
\end{Corollary}

\begin{proof}
Due to Proposition~\ref{beingadomain} the equivalence of (a) and (b) is clear, since a cycle of length $n$ is a complete graph if and only if $n=3$. It also follows from   Proposition~\ref{beingadomain}  that whenever $J_G$ is a prime ideal, then $J_G$ is Cohen--Macaulay, because if each of the components of $G$ is a complete graph, then the binomial edge ideal of each component is the ideal of $2$-minors of a $2\times k$-matrix for some $k$, and these ideals are known to be Cohen--Macaulay. Since $J_G$ is unmixed if $S/I_G$ is Cohen--Macaulay, all implications follow once it is  shown that (c) implies (b). One of the minimal prime ideals of $G$ is $P_{\emptyset}(G)$ and $\dim S/P_{\emptyset}(G)=n+1$. Now let $S\subset [n]$ with $S\neq \emptyset$. We may assume that we have labeled the edges of the cycle counterclockwise, and that
\[
S=\Union_{i=1}^r[a_i,b_i]\quad \text{with}\quad 1=a_1\leq b_1<a_2\leq b_2<\cdots a_r\leq b_r<n.
\]
Then $c(S)=r$, and $\dim S/P_S(G)=n-|S|+c(S)=n-\sum_{i=1}^r(b_i-a_i)-r+r\leq n$. Thus if $J_G$ is unmixed, then $P_{\emptyset}(G)$ is the only minimal prime ideal of $J_G$, and hence since $J_G$ is reduced it follows that $J_G$ is a prime ideal, as required.
\end{proof}

Now let $G$ be an arbitrary simple graph. Which of the ideals $P_S(G)$ are minimal prime ideals of $J_G$? The following result helps to find them.
\begin{Proposition}
\label{minimal}
Let $G$ be a simple graph on $[n]$, and let $S$ and  $T$  be subsets of $[n]$. Let $G_1,\ldots,G_s$  be the connected components of $G_{[n]\setminus S}$, and $H_1,\ldots, H_t$ the connected components of $G_{[n]\setminus T}$.
Then
$P_T(G)\subset P_S(G)$, if and only if  $T\subset S$ and for all $i=1,\ldots,t$ one has $V(H_i)\setminus S\subset V(G_j)$ for some $j$.
\end{Proposition}

\begin{proof}
For a subset $U\subset  [n]$ we let $L_U$ be the ideal generated by the variables $\{x_i,y_i\: \; i\in U\}$. With this notation introduced we have
$P_S(G)=(L_S, J_{\tilde{G_1}},\ldots, J_{\tilde{G}_{s}})$ and $P_T(G)=(L_T, J_{\tilde{H_1}},\ldots, J_{\tilde{H}_{t}})$. Hence  it follows that $P_T(G)\subset P_S(G)$, if and only if $T\subset S$ and $(L_S,  J_{\tilde{H}_1},\ldots, J_{\tilde{H}_{t}})\subset (L_S,J_{\tilde{G}_1},\ldots, J_{\tilde{G}_{s}})$.

Observe that  $(L_S, J_{\tilde{H_1}},\ldots, J_{\tilde{H}_{t}})=(L_S, J_{\tilde{H}_1'},\ldots, J_{\tilde{H}_t'})$ where $H_i'=(H_i)_{[n]\setminus S}$. It follows  that $P_T(G)\subset P_S(G)$ if and only if  $(L_S, J_{\tilde{H}_1'},\ldots, J_{\tilde{H}_t'})\subset (L_S, J_{\tilde{G_1}},\ldots, J_{\tilde{G}_{s}})$ which is the case if and only if $(J_{\tilde{H}_1'},\ldots, J_{\tilde{H}_t'})\subset (J_{\tilde{G_1}},\ldots, J_{\tilde{G}_{s}})$, because the generators of the ideals  $(J_{\tilde{H}_1'},\ldots, J_{\tilde{H}_t'})$ and $(J_{\tilde{G_1}},\ldots, J_{\tilde{G}_{s}})$ have no variables in common with the $x_i$ and $y_i$ for $i\in S$.

Since $V(Hi')=V(H_i)\setminus S$, the assertion will follow once we have shown the following claim: let $A_1, \ldots, A_s$ and $B_1,\ldots,B_t$ be pairwise disjoint  subsets of $[n]$. Then $$(J_{\tilde{A}_1},\ldots, J_{\tilde{A}_s})\subset (J_{\tilde{B}_1},\ldots, J_{\tilde{B}_t}),$$ if and only if for each $i=1,\ldots, s$ there exists a $j$ such that $A_i\subset B_j$.

It is obvious that if the conditions on the $A_i$ and $B_j$ are satisfied, then we have the desired inclusion of the corresponding ideals.

Conversely, suppose that $(J_{\tilde{A}_1},\ldots, J_{\tilde{A}_s})\subset (J_{\tilde{B}_1},\ldots, J_{\tilde{B}_t})$. Without loss of generality we may assume that $\Union_{j=1}^tB_j=[n]$. Consider the surjective $K$-algebra homomorphism
\[
\epsilon\: S\to K[\{x_i,x_iz_1\}_{i\in B_1},\ldots, \{x_i,x_iz_t\}_{i\in B_t}]\subset K[x_1,\ldots,x_n,z_1,\ldots,z_t]
\]
with $\epsilon(x_i) =x_i$ for all $i$ and $\epsilon(y_i)=x_iz_j$  for $i\in B_j$ and $j=1,\ldots,t$. Then
\[
\Ker(\epsilon)=(J_{\tilde{B}_1},\ldots, J_{\tilde{B}_t}).
\]
Now fix one of the sets $A_i$ and let $k\in A_i$. Then $k\in B_j$ for some $k$. We claim that $A_i\subset B_j$. Indeed, let $\ell\in A_i$ with $\ell\neq k$ and suppose that $\ell\in B_r$ with $r\neq j$. Since $x_ky_\ell-x_\ell y_k\in J_{\tilde{A}_i}\subset (J_{\tilde{B}_1},\ldots, J_{\tilde{B}_t})$, it follows that $x_ky_\ell-x_\ell y_k\in \Ker(\epsilon)$, so that $0=\epsilon(x_ky_\ell-x_\ell y_k)=x_kx_\ell z_j-x_kx_\ell z_r$, a contradiction.
\end{proof}

Let $G_1,\ldots,G_r$ be the connect components of $G$.
Once we know the minimal prime ideals of $J_{G_i}$ for each $i$ the minimal prime ideals of $J_G$ are known, Indeed, since the ideals $J_{G_i}$ are ideals in different sets of variables, it follows that the minimal prime ideals of $J_G$ are exactly the ideals $\sum_{i=t}^rP_i$ where each $P_i$ is a minimal prime ideal of $J_{G_i}$.

The next results detects the minimal prime ideals of $J_G$ when $G$ is connected.

\begin{Corollary}
\label{minimalprime}
Let $G$ be a connected simple graph on the vertex set $[n]$, and $S\subset [n]$. Then $P_S(G)$ is a minimal prime ideal of $J_G$ if and only if $S=\emptyset$, or $S\neq\emptyset$ and for each $i\in S$ one has $c(S\setminus \{i\})<c(S)$.
\end{Corollary}

In the terminology of graph theory, the corollary says that if $G$ is a connected graph, then $P_S(G)$ is a minimal prime ideal of $J_G$, if and only if each $i\in S$ is a cut-point of the graph $G_{([n]\setminus S)\union\{i\}}$.

\begin{proof}[Proof of \ref{minimalprime}]
Assume that $P_S(G)$ is a minimal prime ideal of $J_G$. Let $G_1,\ldots,G_r$ be the connected components of $G_{[n]\setminus S}$. We distinguish
several cases.

Suppose that there is no edge $\{i,j\}$ of $G$ such that $j\in G_k$ for some $k$. Set $T=S\setminus \{i\}$. Then the connected components of $G_{[n]\setminus T}$ are $G_1,\ldots,G_r,\{i\}$. Thus $c(T)=c(S)+1$. However this case cannot happen, since  Proposition~\ref{minimal} would imply that $P_T(G)\subset P_S(G)$.

Next suppose that there exists exactly one $G_k$, say $G_1$,  for which there exists $j\in G_1$ such that $\{i,j\}$ is an edge of $G$. Then the connected components of $G_{[n]\setminus T}$  are $G_1', G_2,\ldots, G_r$ where $V(G_1')=V(G_1)\union\{i\}$. Thus $c(T)=c(S)$. Again, this case cannot happen since Proposition~\ref{minimal} would imply that $P_T(G)\subset P_S(G)$.

It remains the case that there are at least two components, say $G_1,\ldots,G_k$, $k\geq 2$, and $j_\ell\in G_\ell$ for $\ell=1,\ldots,k$ such that $\{i,j_\ell\}$ is an edge of $G$. Then the connected components of $G_{[n]\setminus T}$ are $G_1',G_{k+1},\ldots, G_r$, where $V(G_1')=\Union_{\ell=1}^kV(G_\ell)\union\{i\}$. Hence in this case $c(T)<c(S)$.

Conversely, suppose that $c(S\setminus\{i\})<c(S)$ for all $i\in S$. We want to show that $P_S(G)$. Suppose this is not the case. Then there exists a proper subset $T\subset S$ with $P_T(G)\subset P_S(G)$.  We choose $i\in S\setminus T$. By assumption, we have $c(S\setminus\{i\})<c(S)$. The discussion of the three cases above show that we may assume that $G_1',G_{k+1},\ldots, G_r$ are the components of $G_\{[n]\setminus\{i\})$ where $V(G_1')=\Union_{\ell=1}^kV(G_\ell)\union\{i\}$ and where $k\geq 2$. It follows that $G_{[n]\setminus T}$ has one connected component $H$ which contains $G_1'$ Then $V(H)\setminus S$ contains the subsets $V(G_1)$ and $V(G_2)$. Hence $V(H)\setminus S$ is not contained in any $V(G_i)$. According to Proposition~\ref{minimal}, this contradicts the assumption that $P_T(G)\subset P_S(G)$.
\end{proof}

As an example of Corollary~\ref{minimalprime} consider again the cycle $G$ of length $n$. Then,  besides of the prime ideal $P_\emptyset(G)$ which is of height $n-1$, the only other minimal prime ideals are the ideals $P_S(G)$ where $|S|>1$ and and no two elements $i,j\in S$ belong to the same edge of $G$. Each of these prime ideals has height $n$.

\section{CI-Ideals}

Binomial equations and determinantal ideals are of fundamental importance in the theory of conditional independence. In
this final section we will demonstrate the connection between binomial edge ideals and conditional independence (CI)
statements.

We consider a random vector $X=(X_{0},\ldots,X_{N})$ of $N+1$ discrete random variables, where the random variable
$X_{i}$ takes values in the sets $[d_{i}]$ for some positive integers $d_{i}\in\NN$. Then $X$ takes values in
$\mathcal{X} := [d_{0}]\times\dots\times[d_{N}]$. A joint probability distribution of $X$ is a non-negative real valued
function $p:\mathcal{X} \to \RR_{\geq 0}$, such that $\sum_{x\in\mathcal{X}} p(x) = 1$. It can be represented by a real
vector $p = (p_{x_{0},\dots,x_{N}})_{x_{0},\dots,x_{N}}\in\RR^{\mathcal{X}}$, where $p_{x_{0},\ldots,x_{N}}$ stands for
the probability of the event $X_{0}=x_{0},X_{1}=x_{1},\ldots,X_{N}=x_{N}$. In the following we will consider polynomial
equations in these $\prod_{i=0}^{N} d_{i}$ indeterminates, denoting $\CC[p_{x} : x\in \mathcal{X}]$ the ambient
polynomial ring.

For any subset $S\subseteq\{0,\dots,N\}$ we write $X_{S}$ for the collection of random variables $\{X_{i}:i\in S\}$.
Then $X_{S}$ is a random variable on the smaller state space $\mathcal{X}_{S} = \times_{i\in S}[d_{i}]$.  Given
$x_{T}\in \mathcal{X}_{T}$, we denote $\{X_{T} = x_{T}\} := \{y \in \mathcal{X} : y_{i} = x_{i}, \forall i \in T \}$.
The notation $p(X_{T} = x_{T}) := \sum_{x\in\{X_{T}=x_{T}\}}p_{x}$ is common and convenient and may be abbreviated by
$p(x_{T})$, if no confusion can arise.

Let $S$ and $S'$ be two disjoint subsets of $\{0,\dots,N\}$, let $C\subseteq\mathcal{X}$, and fix a joint probability
distribution $p$.  We say that $X_{S}$ is \emph{conditionally independent} of $X_{S'}$ given $C$ (under $p$) iff $p$
satisfies all equations of the form
\begin{equation}
  \label{eq:CIbinomial}
  p(x_{S}^{\phantom{\prime}},x_{S'}^{\phantom{\prime}};C)p(x_{S}',x_{S'}';C) - p(x_{S}^{\phantom{\prime}},x_{S'}';C)p(x_{S}',x_{S'}^{\phantom{\prime}};C) = 0,
\end{equation}
where $x_{S}^{\phantom{\prime}},x_{S}'\in\mathcal{X}_{S}$, $x_{S'}^{\phantom{\prime}},x_{S'}'\in\mathcal{X}_{S'}$, and
\begin{equation}
  \label{eq:pxSxT}
  p(x_{S},x_{S'};C) := p(\{X_{S}=x_{S}\}\cap\{X_{S'}=X_{S'}\}\cap C) = \sum_{\substack{x\in C:\\ x(i) = x_{S}(i)\text{ for }i\in S,\\ x(i) = x_{S'}(i)\text{ for }i\in S'}} p_{x}
\end{equation}
is the probability that $X$ lies in $C$ and agrees with $x_{S}$ on $S$ and with $x_{S'}$ on $S'$. In this case we write
$\ind{X_{S}}{X_{S'}}[C]$.  If $C=\mathcal{X}$, then it is customary to write $\ind{X_{S}}{X_{S'}}$. Let
$T\subseteq\{0,\dots,N\}$ be disjoint from $S$ and $S'$. If $\ind{X_{S}}{X_{S'}}[\{X_{T}=x_{T}\}]$ holds for all
$x_{T}\in\mathcal{X}_{T}$ we write $\ind{X_{S}}{X_{S'}}[X_{T}]$.

An ideal $I$ which is generated by a collection of equations of the form \eqref{eq:CIbinomial} is called a \emph{\CI
  ideal}.  Here, equations \eqref{eq:CIbinomial} are seen as equations among the elementary probabilities $p_{x}$ via
the relations \eqref{eq:pxSxT}.
Note that $I$ is homogeneous. We can identify probability distributions satisfying the equations of $I$ with those
points of the projective variety of $I$ which have real nonnegative homogeneous coordinates.

\begin{Example}
  {\em Consider for a simple example $N=2$ and binary variables $d_{0} = d_{1} = d_{2} = 2$. The polynomial ring is given
    as $\CC[p_{111},p_{112},p_{121},p_{122},p_{211},p_{212},p_{221},p_{222}]$. The conditional independence
    $\ind{X_{0}}{X_{1}}[X_{2}]$ describes the binomial ideal
  \begin{equation*}
    I_{\ind{X_{0}}{X_{1}}[X_{2}]} = \left( p_{111}p_{221}-p_{121}p_{211}, p_{112}p_{222}-p_{122}p_{212}\right)
  \end{equation*}
  In contrast to that, the independence $\ind{X_{0}}{X_{1}}$ is given by the principal ideal
  \begin{equation*}
    I_{\ind{X_{0}}{X_{1}}} = \left( (p_{111}+p_{112})(p_{221}+p_{222}) - (p_{211}+p_{212})(p_{121}+p_{122}) \right).
  \end{equation*}
  }
\end{Example}

\begin{Remark}
{\em   A conditional independence $X_{S} \Perp X_{S'} | C$ is usually defined differently: One requires
  \begin{equation}
    \label{eq:classicalCI}
    p(X_{S}=x_{S},X_{S'}=x_{S'}|X\in C) = p(X_{S}=x_{S}|X\in C) p(X_{S'}=x_{S'}|X\in C)
  \end{equation}
  for all $x_{S}\in\mathcal{X}_{S}$ and $x_{S'}\in\mathcal{X}_{S'}$.  Here,
  \begin{equation*}
    p(X_S=x_{S},X_{S'}=y_{S'}|X\in C) = \frac{p(X_S=x_{S},X_{S'}=y_{S'},X\in C)}{p(X\in C)},
  \end{equation*}
  and so on.  However, equation \eqref{eq:classicalCI} is not well defined if $p(X\in C)$ is zero, while equation
  \eqref{eq:CIbinomial} is defined for all joint distributions $p$.  It is an easy exercise to prove that equations
  \eqref{eq:CIbinomial} and \eqref{eq:classicalCI} are equivalent if $p(X\in C)$ is nonzero.}
\end{Remark}



We will now discuss a special case which makes it possible to apply the results of the first three sections.  Namely, we
assume $d_{0}=2$, i.e., $X_{0}$ is considered to be binary.  In this case we can arrange the elementary probabilities
$p_{x}$ in a $2\times d_{1}\dots d_{N}$-matrix, where the columns are indexed by the state space $\mathcal{X}_{[N]}$ of
$X_{[N]}=(X_{1},\dots,X_{N})$.  The basic observation is that every $2$-minor corresponds to one \CI statement; namely,
the minor
\begin{equation*}
  p_{1x} p_{2x'} - p_{2x} p_{1x'}
\end{equation*}
of the two columns corresponding to $x,x'\in\mathcal{X}_{[N]}$ expresses exactly the \CI statement
\begin{equation*}
  \ind{X_{0}}{X_{[N]}}[\left\{X_{[N]}\in\{x,x'\}\right\}].
\end{equation*}
In this way we can associate a collection of \CI statements to every graph on the vertex set $\mathcal{X}_{[N]}$.

Until now we did not use of the fact that $X_{[N]}$ is a product of several random variables.  Now let $S\cup T$ be a
(disjoint) partition of $[N]$ and consider the \CI statement
\begin{equation}
  \label{eq:X0indST}
  \ind{X_{0}}{X_{S}}[X_{T}].
\end{equation}
For simplicity we assume that $S=\{1,\dots,s\}$ for a moment.  Then \eqref{eq:X0indST} is equivalent to the equations
\begin{equation*}
  p_{1x_{S}x_{T}} p_{2x_{S}'x_{T}} - p_{1x_{S}'x_{T}} p_{2x_{S}x_{T}}
\end{equation*}
for all $x_{S},x_{S}'\in\mathcal{X}_{S}$ and $x_{T}\in\mathcal{X}_{T}$.  These equations come from all $2$-minors with
columns $x,x'\in\mathcal{X}_{[N]}$ such that $x$ and $x'$ agree on their $T$-components.  This means that we can
associate with \eqref{eq:X0indST} the graph on $\mathcal{X}_{[N]}$ with edges
\begin{equation*}
  E(G) = \{(x,x') : \text{ $x,x'\in\mathcal{X}_{[N]}$ agree on $T$}\}.
\end{equation*}
More generally, when we have a collection $\mathcal{C} = \{\ind{X_{0}}{X_{S_{i}}}[X_{T}]\}$ of \CI statements
corresponding to disjoint partitions $S_{i}\cup T_{i}$ of $[N]$, we can associate a graph $G_{i}$ with every single
statement.  If we define a graph $G$ on $\mathcal{X}_{[N]}$ by $E(G) = \Union_{i}E(G_{i})$, then the binomial edge ideal
of $G$ equals the \CI ideal of $\mathcal{C}$.

\CI statements of the form under consideration have the following natural interpretation in probabilistic modeling: We
consider $X_{0}$ as the output node of a system which receives input from $X_{1},\dots,X_{N}$. Then we can ask how much
information is lost when certain input nodes are not available. If $\ind{X_{0}}{X_{S'}}[X_{T}]$, then all the relevant
information can be reconstructed from $X_{T}$ alone: The system can dispense with the information from $X_{S'}$.  In
this way, a collection of \CI statements can be used to model a notion of robustness of probabilistic computation
\cite{AK}.  Because of this interpretation we introduce the following notation:
\begin{Definition}
  A collection of \CI statements induced as above by a set of disjoint partitions $S_{i}\cup T_{i} = [N]$ will be called
  a \emph{robustness specification}.
\end{Definition}
Theorems~\ref{radical} and~\ref{intersection} imply:
\begin{Corollary}
  The \CI ideal of a robustness specification 
  with binary output is a radical ideal.
\end{Corollary}

Now fix a robustness specification $\mathcal{C}$.  Owing to Theorem~\ref{intersection}, each minimal prime is given by a
subset $S\subseteq\mathcal{X}_{[N]}$ which satisfies the conditions of Corollary \ref{minimalprime}.  Such a subset $S$
defines events with zero probability: $p(X_{[N]}\in S) = 0$ if $p\in V(P_{S}(G))$, where $G=G_{\mathcal{C}}$.  In the
language of statistical modeling, $S$ is a set of structural zeros.
\begin{Corollary}
  Let $I$ be the \CI ideal of a robustness specification.  Each minimal prime $P$ of $I$ is characterized by a set $S$
  of structural zeros in the distribution of $X_{[N]}$ which is common to all probability distributions lying in the
  component corresponding to $P$.  The possible sets $S$ are characterized by Corollary \ref{minimalprime}.
\end{Corollary}

The binomial generators $J_{\tilde{G}_{1}},\dots,J_{\tilde{G}_{c(S)}}$ in $P_{S}(G)$ also have a nice statistical
interpretation: Namely $J_{\tilde{G}_{i}}$ expresses the \CI statement
\begin{equation*}
  \ind{X_{0}}{X_{[N]}}[\left(X_{[N]}\in G_{i}\right)].
\end{equation*}
This means: If we know $S$, then the knowledge in which component of $G_{[N]\setminus S}$ the random vector $X_{[N]}$
lies contains all the relevant information about $X_{0}$.  Once we know this component, the conditional probability
distribution of $X_{0}$ is independent of any further information we may obtain.  In other words, if we know $G$ and
$S$, then we can define a random variable $C$ which maps every outcome of $X$ with nonzero probability to the
corresponding component in $[c(S)]$.  We then have $\ind{X_{0}}{X_{[N]}}[C]$, a fact which can be depicted by the
following Markov chain
\begin{equation*}
  X_{[N]} \longrightarrow C \longrightarrow X_{0}.
\end{equation*}

This corresponds to the classical result that each irreducible component of a binomial ideal is essentially a toric
variety \cite{ES}, and in particular each irreducible component has a rational parametrization.  The most natural such
parametrization in the statistical setting is the following: $p$ factors as a product of a distribution on the connected
components $G_{1},\ldots,G_{c(S)}$ and a distribution of $X_{0}$ for each of the connected components.  This should be
compared to the dimension $n - |S| + c(S)$ in Lemma~\ref{dimension}.

Each binomial ideal $I \subset \CC[p_{x} : x\in \mathcal{X}]$ has the toric ideal $I :
(\prod_{x\in\mathcal{X}}p_{x})^{\infty}$ as a minimal prime. It corresponds to $S=\emptyset$, and all distributions with
full support ($p(x)>0$ for all $x\in \mathcal{X}$) satisfying the robustness specification are contained in the toric
variety. We obtain the following
\begin{Corollary}
  Let $p$ be a probability distribution satisfying the robustness specification $\mathcal{C} =
  \{\ind{X_{0}}{X_{S_{i}}}[X_{T_{i}}]: i=1,\dots,r\}$.  If $p$ has full support (i.e., $p_{x}>0$ for all
  $x\in\mathcal{X}$), then
  \begin{equation*}
    \ind{X_{0}}{X_{\cup_{i}S_{i}}}[X_{\cap_{i}T_{i}}].
  \end{equation*}
  In particular, if $\cup_{i}S_{i} = [N]$ then $\ind{X_{0}}{X_{[N]}}$ and $X_{0}$ is unconditionally independent of the
  input.
\end{Corollary}
\begin{Remark}
{\em   It is easy to prove this corollary directly using the \emph{intersection axiom} \cite{DSS}.}
\end{Remark}

This result is not surprising: If any combination of inputs in $\mathcal{X}_{[N]}$ is possible, then we can't deduce any
missing information. Any distribution where $X_{0}$ is robust against perturbation of the inputs must make use of
features of the input statistics.

\begin{Examples}
  {\em Fix $k\in[N]$ and consider the collection of \CI statements
    \begin{equation}
     \label{eq:krobuststatements}
     \left\{ \ind{X_{0}}{X_{S}}[X_{T}] : S \in \binom{[N]}{k} \right\}
  \end{equation}
  induced by all $k$-element subsets of $[N]$. Consider the graph $G_{k}$ with vertices $\mathcal{X}_{[N]}$ and edges
  between any $x$ and $y$ which differ in at most $k$ components. In other words, $\{x,y\}\in E(G_{k})$ if and only if
  the \emph{Hamming distance} between $x$ and $y$ is at most $k$. The \CI ideal for the
  statements \eqref{eq:krobuststatements} is the binomial edge ideal of $G_{k}$.

  \noindent
  (a) If $k=1$ and $d_{i}=2$, for all $i\in [N]$ we find the graph of the $N$-cube.

  \noindent
  (b) If $k=1$ and $N=2$ we have just two \CI statements:
  \begin{equation*}
    \ind{X_{0}}{X_{1}}[X_{2}]\text{ and }\ind{X_{0}}{X_{2}}[X_{1}].
  \end{equation*}
  These statements have been studied by A.~Fink \cite{AF}. In this case the minimal primes can be seen to correspond to
  bipartite graphs $\Gamma$ such that every connected component is a complete bipartite graph. The two groups of
  vertices in these graphs are $[d_{1}]$ and $[d_{2}]$.  The corresponding prime is minimal if each vertex belongs to at
  least one edge. Such bipartite graphs are in bijection with pairs of partitions $[d_{1}] = I_{1}\cup\dots\cup I_{c}$
  and $[d_{2}] = J_{1}\cup\dots\cup J_{c}$, where $c$ is the number of connected components of $\Gamma$, and $I_{i}$
  resp. $J_{i}$ are the vertices in the $i$th component of $\Gamma$.  Then $S = \mathcal{X}_{[N]}\setminus
  \cup_{i=1}^{c}(I_{i}\times J_{i})$ gives the link with our notation.  In other words, the vertices of the connected
  components $G_{1},\dots,G_{c(S)}$ are given by $V(G_{i}) = I_{i}\times J_{i}$.

  \noindent
  (c) The considerations of (b) generalize to the case $k=N-1$: As above, the minimal primes correspond to partitions
  $[d_{i}] = I_{i,1}\cup\dots\cup I_{i,c}$, where $S = \mathcal{X}_{[N]}\setminus
  \cup_{j=1}^{c}(I_{1,j}\times\dots\times I_{N,j})$, and the components of $G_{T}$ satisfy $V(G_{i}) =
  I_{1,j}\times\dots\times I_{N,j}$. We leave the verification of these results as an exercise to the reader.
  Unfortunately, the nice form of the connected components of $G_{T}$ does not generalize for $k < N-1$.}
\end{Examples}

\section*{Acknowledgment}
\label{sec:acknowledgment}
We wish to thank Seth Sullivant who established the contact between the two groups of authors. The last two authors
thank Nihat Ay for pointing them at conditional independence problems in robustness theory.

{}

\end{document}